\documentclass{ifacconf}

\usepackage{xcolor}
\usepackage{graphicx}
\usepackage{xspace}
\usepackage{algorithm}
\usepackage{algpseudocode}
\algrenewcommand\algorithmicrequire{\textbf{Input:}}

\usepackage{amsmath,amsfonts,amssymb}
\usepackage{bm,fixmath} 
\usepackage{tikz}
\usetikzlibrary{arrows, positioning, calc}

\newtheorem{assumption}{Assumption}

\DeclareMathOperator*{\argmin}{arg\,min}

\newcommand{\N}{\mathbb{N}}
\newcommand{\R}{\mathbb{R}}

\newcommand{\norm}[1]{\left\lVert#1\right\rVert}

\newcommand{\x}{\mathbold{x}}

\newcommand{\w}{\mathbold{w}}

\newcommand{\A}{\mathbold{A}}
\newcommand{\B}{\mathbold{B}}

\newcommand{\G}{\mathbold{G}}

\usepackage{accents}
\newcommand{\lmax}{\bar{\lambda}}
\newcommand{\lmin}{\underline{\lambda}}

\usepackage{natbib} 
\begin{document}
\begin{frontmatter}

\title{Adaptive Online Optimization for Microgrids with Renewable Energy Sources\thanksref{footnoteinfo}} 

\thanks[footnoteinfo]{
Corresponding author: N. Bastianello (\texttt{nicolba@kth.se})\\
The work of W.J.A. van Weerelt was supported by the Digital Futures Summer Research Internship Program.\\
The work of A. Fontan was supported by the Knut and Alice Wallenberg Foundation Wallenberg AI, Autonomous Systems and Software Program (WASP) funded by the
Knut and Alice Wallenberg Foundation.\\
The work of N. Bastianello was partially supported by the EU Horizon Research and Innovation Actions program under Grant 101070162.}

\author[First]{Wouter J.A. van Weerelt} 
\author[Second]{Angela Fontan} 
\author[Second]{Nicola Bastianello}

\address[First]{Department of Mathematics, KTH Royal Institute of Technology, Stockholm, Sweden}
\address[Second]{School of Electrical Engineering and Computer Science, and Digital Futures, KTH Royal Institute of Technology, Stockholm, Sweden}

\begin{abstract} 
In this paper we propose a novel adaptive online optimization algorithm tailored to the management of microgrids with high renewable energy penetration, which can be formulated as a constrained, online optimization problem.
The proposed algorithm is characterized by a control-based design that applies the internal model principle, and a system identification routine tasked with identifying such internal model. In addition, in order to ensure the constraints are verified, we integrate a projection onto the constraint set.
We showcase promising numerical results for the microgrid use case, highlighting in particular the enhanced adaptability of the proposed algorithm to changes in the internal model. The performance of the proposed algorithm is shown to outperform state-of-the-art alternative in the long-term, ensuring efficient management of the grid.
\end{abstract}


\end{frontmatter}
\section{Introduction}\label{sec:introduction}
The power grid is undergoing significant changes, driven by the increasing penetration of renewable energy sources and the need to mitigate the climate crisis \citep{OO_survey}. In fact, between 2012 and 2023, renewable electricity expanded at a compound rate of 5.9\% compared to the 1.3\% of non-renewables \citep{IRENA}. Likewise, renewable energy microgrids have been shown to improve livelihoods, contribute to economic growth as well as enhance both food security and health \citep{microgrid-survey}.
Efficiently exploiting renewable sources to satisfy user demand, while reducing the use of fossil fuels as back-up sources, is however a complex problem.
A successful approach to addressing this objective is to leverage \textit{online optimization} \citep{OPF_dallanese,OO_survey}. The key idea is to formulate the grid management problem as an online problem, whose time-variability models the (possibly high frequency) fluctuations in user demand and renewables availability.

Different power grid management problems have been formulated as online optimization problems, and suitable online algorithms have been proposed to solve them.
For example, typical optimal power flow problems, which can be solved with online gradient algorithms \citep{OPF_dallanese, radialOPF} or ADMM algorithms \citep{ADMM}. Likewise, multi-period optimization problems can be solved with other algorithms ranging from dynamic, distributed methods such as \citep{lesage} to algorithms grounded in game theory \citep{DR_survey}.

These results showcase the promise of the online optimization framework for power grid management. However, there is potential for performance improvements by designing tailored online algorithms.
Indeed, online algorithms can be largely divided into two groups: \textit{unstructured}, which are ``model-agnostic" in that they do not exploit information on the time-variability of the problem, and \textit{structured}, which are instead ``model-based" \citep{simonetto}. By their design, structured methods can achieve higher performance than the unstructured algorithms currently being applied to grid management problems.
Therefore the goal of this paper is to showcase the use of a novel, structured algorithm to the solution of these problems.

Our proposed algorithm in particular is designed adopting a ``control-based" approach, which has been successfully leveraged to achieve near-zero asymptotic optimality error in a variety of settings, including unconstrained \citep{Bastianello2024}, constrained \citep{casti2023}, distributed \citep{vanweerelt2025} and stochastic \citep{Casti2024}.
The success of these algorithms comes from exploiting a model of the online problem to enhance performance. However, the knowledge of this internal model might not be available, especially in a power grid context.
Therefore, in this paper we build on the algorithm proposed in \citep{weerelt_selfidentifying_2025} which adapts to the online problem by concurrently identifying the internal model, and applying a control-based algorithm using the identified model.
We focus in particular on a demand response application in the context of microgrids with renewable energy sources.
Our contributions are summarized as follows:
\begin{itemize}
    \item We propose a novel, adaptive online algorithm designed using control theory. The algorithm, differently from previous works, incorporates an identification routine which adapts online to changes in the demand and availability profiles.

    \item We evaluate the performance of the proposed algorithm for a microgrid which integrates renewable energy sources. The results show how our algorithm outperforms, in the long-term, unstructured methods.
\end{itemize}

\section{Problem Formulation}\label{sec:problem}

In this section we formulate the specific power grid management problem that we address. We focus on deploying online optimization algorithms in the context of demand response in power distribution systems \citep{lesage}. More precisely, we consider the management of controllable distributed energy resources (DERs) on a microgrid, as done in \citet{Ospina}.

A microgrid can be seen as a small portion of a larger power distribution network \citep{Bolognani}. Such a microgrid hosts a number of loads and several small modular generation systems, which can provide power and heat \citep{Lopes}. Examples of microgrids include small urban areas, a shopping center or an industrial park. With the increased penetration of renewables in microgrids, more controllable DERs must be used to counteract the volatility and uncertainty of the grid \citep{OO_survey}. Examples of DERs in a residential microgrid may include solar panels and battery energy storage systems. Typically, microgrids are considered either directed or undirected graphs with certain topology, connected to the transmission grid at a point of common coupling, or PCC.

We translate the problem of managing the microgrid into an online optimization problem as follows.
We consider $n$ controllable inputs $\boldsymbol x \in \mathbb{R}^n$ and unknown exogenous inputs $\boldsymbol w \in \mathbb{R}^w$. As in \citep{Ospina}, we model the microgrid as an algebraic map $\boldsymbol y=\mathcal{M}(\boldsymbol x, \boldsymbol w)$ where $\mathcal{M} : \mathbb{R}^n \times \mathbb{R}^w \rightarrow \mathbb{R}^y $ is well defined. 
The problem then can be expressed as:
\begin{equation}\label{eq:objective}
    \boldsymbol x_{*,k} \in \min_{\boldsymbol x\in \mathcal{X}_k}:=C_k(\mathcal{M}(\boldsymbol x,\boldsymbol w_k)+\boldsymbol U_k(\boldsymbol x).
\end{equation}
Here, we consider $k\in \mathbb{N}$ to be the time index, $\mathcal{X}_k \subseteq \mathbb{R}^n$ to be a time-varying, convex constraint set for inputs $\boldsymbol{x}$, $\boldsymbol x \mapsto U_k(\boldsymbol x)$ to be a cost associated with inputs, and $\boldsymbol y \mapsto C_k(\boldsymbol y)$ to be a cost associated with outputs. The time-varying constraint set depends on physical constraints of the system and the DERs.
We define the algebraic map as $\mathcal{M}(\boldsymbol x, \boldsymbol w_k) =\boldsymbol J \boldsymbol x+\boldsymbol H\boldsymbol w_k$, which is built based on a linearization of the power flow equations \citep{Ospina,OPF_dallanese}. We remark that the matrices $\boldsymbol J$ and $\boldsymbol H$ encode physical information of microgrid and are dependent on the grid's components and topology. Furthermore, $\boldsymbol y_k$ represents the net real power exchange at some points of common coupling (PCCs). It is further assumed that the power consumed by the uncontrollable loads $\boldsymbol{w}$ cannot be individually measured, meaning that measurements of $\boldsymbol y_k$ are instead available from meters and sensing units. The function $\boldsymbol U_k(\boldsymbol x)$ in \eqref{eq:objective} represents user dissatisfaction which is considered to be strongly-convex and quadratic. Finally, we assume the cost associated with the outputs to be defined as: $C_k(\boldsymbol y_k) = \frac{\beta}{2}\norm{\boldsymbol y_k-\boldsymbol y_{k,\text{ref}}}^2$, where $\boldsymbol y_{k,\text{ref}}$ is a time-varying demand response setpoint for the PCCs and $\beta>0$ is a given parameter.

With these assumptions in place, problem~\eqref{eq:objective} becomes the following:
\begin{equation}
    \boldsymbol x_{*,k} \in \min_{\boldsymbol x\in \mathcal{X}_k}:=\frac{\beta}{2}\norm{\boldsymbol J \boldsymbol x+\boldsymbol H\boldsymbol w_k-\boldsymbol y_{k,\text{ref}}}^2+\boldsymbol U_k(\boldsymbol x).
\end{equation}
However, since we assume uncontrollable loads to be unmeasurable, we need to reformulate this problem into:

\begin{equation}\label{eq:true_objective}
        \boldsymbol x_{*,k} \in \min_{\boldsymbol x\in \mathcal{X}_k}:=\frac{\beta}{2}\norm{\boldsymbol y_k -\boldsymbol y_{k,\text{ref}}}^2+\boldsymbol U_k(\boldsymbol x).
\end{equation}

\smallskip

Problem~\eqref{eq:true_objective} is a constrained, online optimization problem, with a smooth cost function and a convex constraint set.
This problem in principle can be solved with an online version of projected gradient descent \citep{Nocedal2006}. However, the resulting algorithm would be unstructured, and in the next section we design a structured algorithm instead, aiming for improved performance.

\section{Proposed Algorithm}\label{sec:algorithm}
As discussed in section~\ref{sec:introduction}, in this paper we propose a control-based algorithm design to solve~\eqref{eq:true_objective}. The problem has a quadratic cost but also a constraint set. Therefore, the control-based algorithm proposed in \citep{Bastianello2024} (and the counterpart for linearly constrained problems of \citep{casti2023}) does not apply to it directly.
In section~\ref{subsec:algorithm} thus we revise the algorithm design of \citep{Bastianello2024}, based on the internal model principle, to incorporate a projection onto the constraint set $\mathcal{X}_t$.
Additionally, to overcome the unmeasurability of some loads, in section~\ref{subsec:identification} we design an identification routine for the internal model.
The result is an online algorithm which automatically adapts to the demand and availability profile in the microgrid, while ensuring improved performance through the control-based design.

\subsection{Projected control-based algorithm}\label{subsec:algorithm}
Abstracting from~\eqref{eq:true_objective}, we are interested in solving the online constrained, quadratic problem
\begin{equation}\label{eq:general-problem}
    \x _{*,k} = \argmin_{\x \in \mathcal{X}_k} f_k(\x):=\frac{1}{2}\x^\intercal\A\x +\x^\intercal \mathbold b_k, \quad k \in \mathbb{N}
\end{equation}
under the following assumptions.

\begin{assumption}[Cost and constraint set]\label{as:cost-constraint}
The symmetric matrix $\A$ is such that $\lmin\mathbold{I}\preceq\A =\A^\intercal\preceq \lmax\mathbold I$, with $0<\lmin<\lmax<\infty$. $\mathcal{X}_k \subset \R^n$ is a convex, non-empty set.
\end{assumption}

\begin{assumption}[Model of $\mathbold b_k$]\label{as:model-b}
There exists a rational $\mathcal{Z}$-transform for the sequence $\{\mathbold b_k\}_{k\in\mathbb{N}}$, namely:
\begin{align} \label{eq:Z-transform}
    \mathcal{Z}[\mathbold b_k] = \B(z) = \frac{\B_N(z)}{B_D(z)}, \hspace{0.2cm} B_D(z) = z^m+\sum_{i=0}^{m-1}d_iz^i,
\end{align}
with $\B_N(z) = \sum_{i=0}^p \mathbold u_iz^i$, $p \leq m$. We assume that the poles of $B_D(z)$ are all marginally or asymptotically stable.
\end{assumption}

Assumption~\ref{as:cost-constraint} ensures that the cost functions $\{f_k\}_{k\in\mathbb{N}}$ are $\lmin$-strongly convex and $\lmax$-smooth for any time $k\in\mathbb{N}$, and that a (unique) solution to~\eqref{eq:general-problem} always exists. This assumption is verified in the setting of~\eqref{eq:true_objective}.
Assumption~\ref{as:model-b} ensures that the sequence of linear terms $\{ \boldsymbol{b}_k \}_{k \in \N}$ admits a model which can serve as internal model when designing the control-based algorithm.

Disregarding the constraints in~\eqref{eq:general-problem} for now (\textit{i.e.} $\mathcal{X}_k \equiv \R^n$), we can apply \citep{Bastianello2024} to design the following algorithm:
\begin{subequations}\label{eq:baseline-algorithm}
\begin{align}
    \boldsymbol w_{k+1} &= \left( \boldsymbol F \otimes \boldsymbol I \right) \boldsymbol w_k + \left( \boldsymbol G \otimes \boldsymbol I \right) \nabla f_k (\boldsymbol x_k) \\
    \boldsymbol x_{k+1} &= \left( \boldsymbol K \otimes \boldsymbol I \right)\boldsymbol w_{k+1}
\end{align}
\end{subequations}
where $\boldsymbol{w}$ is the internal state of the algorithm, and
\begin{equation}\label{eq:state-space-realization}
\begin{split}
	\boldsymbol{F} &= \begin{bmatrix}
		0 & 1 & 0 & \cdots \\
		& & \ddots & \\
		0 & \cdots & 0 & 1 \\
		-d_0 & \cdots & \cdots & -d_{m-1}
	\end{bmatrix},\quad
	\G= \begin{bmatrix}
		0\\
		\vdots\\
		0\\
		1
	\end{bmatrix}\\
    \boldsymbol{K} &=\begin{bmatrix}c_0 & c_1 & \cdots & \cdots & c_{m-1}\end{bmatrix}.
\end{split}
\end{equation}
The controller $\mathbold{K}$ can be computed in order to stabilize the closed loop if the internal model $B_D(z)$ is known.

\begin{figure}[!ht]
    \centering
    \includegraphics[width=0.5\linewidth]{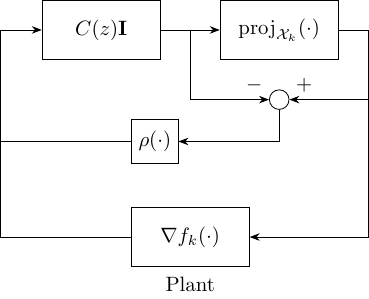}
\caption{Block diagram of algorithm \eqref{eq:proj_CB}}
\label{fig:alg9}
\end{figure}

In the use case described in section~\ref{sec:problem}, the problem is constrained by $\mathcal{X}_k$, and thus~\eqref{eq:baseline-algorithm} does not directly apply.
The idea then is to incorporate into the algorithm a projection onto $\mathcal{X}_k$, yielding the following:
\begin{subequations}\label{eq:projected-algorithm}
\begin{align}
    \boldsymbol w_{k+1} &= \left( \boldsymbol F \otimes \boldsymbol I \right) \boldsymbol w_k + \left( \boldsymbol G \otimes \boldsymbol I \right) \nabla f_k (\boldsymbol x'_k) \\
    \boldsymbol x_{k+1} &= \left( \boldsymbol K \otimes \boldsymbol I \right)\boldsymbol w_{k+1} \\
    \boldsymbol x'_{k+1} &= \operatorname{proj}_{\mathcal{X}_k}(\boldsymbol x_{k+1})
\end{align}
\end{subequations}
where the output of~\eqref{eq:baseline-algorithm} is filtered by the projection, ensuring that $\x'_k$ belongs to the constraint set.
As observed in \citep{casti2023}, however, projections act as nonlinearities similar to saturations. But since~\eqref{eq:baseline-algorithm} is designed using linear control tools, a saturation might significantly impact performance.
Therefore, following the solution proposed in \citep{casti2023}, we add an anti-windup module to~\eqref{eq:projected-algorithm} to reduce the impact of the nonlinear projection. The result is the algorithm:
\begin{subequations}\label{eq:proj_CB}
\begin{align}
    \boldsymbol w_{k+1} &= \left( \boldsymbol F \otimes \boldsymbol I \right) \boldsymbol w_k + \left( \boldsymbol G \otimes \boldsymbol I \right) (\nabla f_k (\boldsymbol x'_k)-\rho(\boldsymbol x'_k-\boldsymbol x_k)) \\
    \boldsymbol x_{k+1} &= \left( \boldsymbol K \otimes \boldsymbol I \right)\boldsymbol w_{k+1}\\
    \boldsymbol x'_{k+1} &= \operatorname{proj}_{\mathcal{X}_k}(\boldsymbol x_{k+1})    
    \end{align}
\end{subequations}
where $\rho > 0$ is the tunable anti-windup parameter. Figure~\ref{fig:alg9} depicts this algorithm as a block diagram.

\smallskip

Algorithm~\eqref{eq:proj_CB} in principle could be applied to solve the microgrid management problem~\eqref{eq:true_objective}. However, Assumption~\ref{as:model-b} does not hold in practice, due to the presence of unmeasurable loads.
This implies that the model of the linear term in~\eqref{eq:true_objective}, which is
$$
    \boldsymbol{b}_k = \beta\boldsymbol J^\intercal(\boldsymbol H\boldsymbol w_k-\boldsymbol y_{k,\text{ref}}) + \boldsymbol{u}_{2,k},
$$
is not known a priori, where $\boldsymbol{u}_{2,k}$ is the linear time-varying term of $U_k$ (See section~\ref{sec:numerics} for the details).
Indeed, $\boldsymbol{b}_k$ depends on the uncontrollable and unmeasurable loads $\w_k$, and on $\boldsymbol J$, $\boldsymbol H$, which encode physical information on the microgrid and might not be available.

\smallskip

The solution then, outlined in the next section, is to gather information on $\boldsymbol{b}_k$ in an online fashion, reconstructing the internal model which can then be used in~\eqref{eq:proj_CB}.

\subsection{Online identification routine}\label{subsec:identification}
\begin{figure}[!ht]
    \centering
    \includegraphics[width=\linewidth]{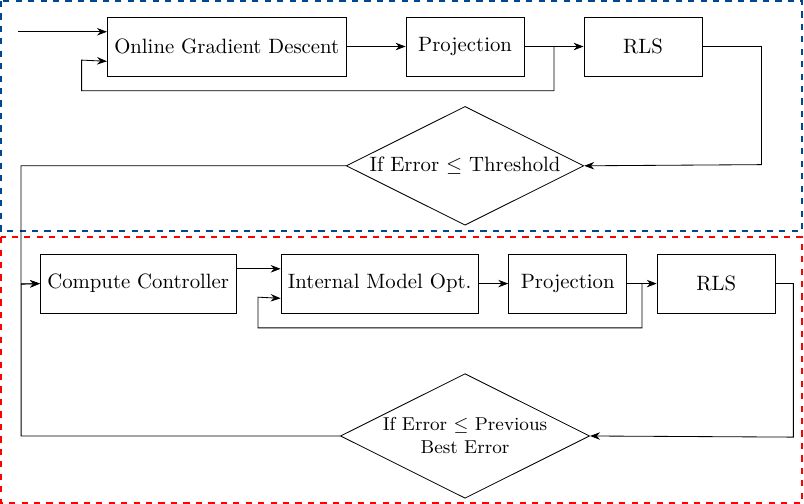}
\caption{Flowchart of the proposed algorithm}
\label{fig:flowchart}
\end{figure}

The goal of this section is to design an online identification routine that reconstructs the internal model of $\boldsymbol{b}_k$.
We start by introducing the approach proposed in \citep{weerelt_selfidentifying_2025}, and then show how to integrate projections into it.
The algorithm of \citep{weerelt_selfidentifying_2025} is split in two phases, depicted in Figure~\ref{fig:flowchart}.
The first phase (top box) involves identifying the internal model based on the output of the online gradient descent (ignoring the projection for now):
\begin{equation}\label{eq:proj_GD}
    \x_{k+1} = \x_k - h \nabla f_k(x_k).
\end{equation}
This can be considered a warm-up phase, where a simple (unstructured) algorithm is applied to start collecting information on the internal model. The specific identification procedure we use is \textit{recursive least squares} (RLS), which is computationally efficient and lends itself to an online implementation \citep{TZLai}.
Once the solution output by the online gradient no longer improves, in the sense that the \textit{a posteriori} error given by RL, \textit{i.e.,} how well the identified internal model coefficients match the true ones, settles on a roughly constant value, the second phase (bottom box in Figure~\ref{fig:flowchart}) is triggered. In this phase, the internal model identified during the first phase is used to computed the controller ($\boldsymbol{K}$) which defines algorithm~\eqref{eq:baseline-algorithm} (again, ignoring the projection for now).
However, the internal model at this stage might not be correctly identified or, more importantly, it might be changing over time. Therefore, in the second phase the recursive least squares identification runs concurrently with the online algorithm, to improve identification accuracy or adapt to changes. The re-identified model is then used to recompute the controller when the identification error settles on an approximately constant value.

We are now ready to concretely outline the identification routine we apply. To apply recursive least squares (RLS), we need a recurrence relating the output of the online algorithm (be it online gradient descent in phase one or the control-based method in phase two) to the internal model coefficients \citep{Guo1994}.
To this end, consider the $\mathcal{Z}$-transform of the output of the online gradient descent in phase one, $\mathcal{Z}\{\x_{k+1}\} = \mathcal{Z}\{\mathbold x_{k} - h \nabla f_k(\x_k) \}$, which yields:
\begin{equation}\label{eq:timeshiftz}
    z\mathbold X(z) = (\mathbold I - h \mathbold A)\mathbold X(z) - h \mathbold B(z).
\end{equation}
Rewriting the $\mathcal{Z}$-transforms as infinite sums \citep{Graf2004} and using Assumption~\ref{as:model-b} yields:
\begin{align*}
    &\sum_{k=0}^\infty \left( \mathbold x_{k+m}  +\sum_{i=0}^{m-1}d_i\mathbold x_{k+i}\right)z^{-k} = \nonumber \\ &\hspace{2.5cm}-h \sum_{k=0}^\infty(\mathbold I-h \mathbold A)^kz^{-k-1} \sum_{j=0}^p \mathbold u_jz^j
\end{align*}
which we simplify to

\begin{equation}
    \mathbold x_{k+m}  +\sum_{i=0}^{m-1}d_i\mathbold x_{k+i} = -h \sum_{j=0}^p(\mathbold I - h \mathbold A)^{k-1+j}\mathbold{u}_jz^j.
\end{equation}
By Assumption~\ref{as:cost-constraint}, $\mathbold A$ is symmetric positive definite, and selecting the stepsize $h < 2 / \lmax$, it holds that:
\begin{align}
    \lim_{k\to\infty}\left(\mathbold I - h \mathbold A\right)^{k-1+j} = \boldsymbol 0.
\end{align}
This results in the following recurrence relation: \begin{align} \label{eq:recurrugly2}
        \mathbold x_{k+m} + \sum_{i=0}^{m-1}d_i\mathbold x_{k+i} = \boldsymbol 0, \hspace{2ex} \forall k \geq m +1,
\end{align}
where recall that $d_i$ are the coefficients of the internal model $B_D(z)$ (see Assumption~\ref{as:model-b}), and $m$ is its order.
It is then possible to apply RLS to this recurrence and identify the coefficients $d_i$'s. In particular, let $\mathbold{d} = [d_0, \ldots, d_{m-1}]^\intercal \in \R^m$ be the vector of coefficients of the polynomial $B_D(z)$. Then RLS identifies it using the recursion:
\begin{equation}\label{eq:RLS_basic}
    \hat{\mathbold {d}}_{k+1} = \hat{\mathbold{d}}_{k} + \mathbold L_k(\mathbold y_k-\mathbold{\phi}_k^\intercal\hat{\mathbold{d}}_k),
\end{equation} 
where $\mathbold y_k$ denotes the observation data, $\mathbold \phi_k$ denotes the regressor, and $\mathbold L_k$ denotes the gain \citep{Guo1994}:
    with $\mathbold L_k$ given by: \begin{align}\label{eq:gain-rls}
    \mathbold L_k = \mathbold P_k\mathbold{\phi}_k (\alpha\mathbold I +\mathbold{\phi}_k^\top \mathbold P_k \mathbold{\phi}_k)^{-1},
\end{align}
and: \begin{align}\mathbold P_{k+1}=\frac{1}{\alpha}\left(\mathbold P_k -\mathbold P_k\mathbold{\phi}_k (\alpha\mathbold I +\mathbold{\phi}_k^\top \mathbold P_k \mathbold{\phi}_k)^{-1}\mathbold{\phi}_k^\top\mathbold P_k\right).\end{align} 
In these equations, $\mathbold P_0 > 0$, and $\alpha \in (0,1)$ is a forgetting factor which weights recent entries more heavily \citep{Guo1994}.

A similar identification routine can be applied during phase two, when the control-based algorithm~\eqref{eq:baseline-algorithm} is being applied.
The controller is designed using the inexact internal model $\hat{B}_D(z) = z^m+\sum_{i=0}^{m-1} \hat{d}_iz^i$, computed during phase one, to yield
\begin{equation}\label{eq:inexact-ctrl}
    C(z) = \frac{C_N(z)}{\hat{B}_D(z)}.
\end{equation}
The resulting closed-loop system is then characterized by the $\mathcal{Z}$-transform:
\begin{equation*}
    \mathbold X(z) = \left(\hat{B}_D(z) \mathbold I-C_N(z)\mathbold A\right)^{-1}\frac{\mathbold B_N(z)C_N(z)}{B_D(z)}
\end{equation*}
and hence
\begin{equation}\label{eq:CB-rec-der}
    B_D(z)\mathbold X(z)  = \left(\hat{B}_D(z) \mathbold I-C_N(z)\mathbold A\right)^{-1}\mathbold B_N(z)C_N(z).
\end{equation}
Using again the infinite sum expression for $\mathcal{Z}$-transforms allows to rewrite~\eqref{eq:CB-rec-der} as: \begin{equation}
\begin{split}
\sum_{k=0}^\infty \Big(\mathbold x_{k+m}
  +\sum_{i=0}^{m-1} d_i \mathbold x_{k+i}\Big) z^{-k} = \\ \qquad \Big(\hat{B}_D(z)\mathbold I - C_N(z)\mathbold A\Big)^{-1} \mathbold B_N(z)\,C_N(z)
\end{split}
\end{equation}
Since the right-hand side of the equation is made up of two anti-causal signals, and $C_N(z)$ is chosen so that $\left(\hat{B}_D(z) \mathbold I-C_N(z)\mathbold A\right)$ is stable, then the resulting recurrence is:
\begin{subequations}
\begin{align}\label{eq:p2_RLS_rec}
           \mathbold x_{k+m}  +\sum_{i=0}^{m-1}d_i\mathbold x_{k+i} = \boldsymbol 0, \hspace{2ex} k\geq m+1.
\end{align}
\end{subequations}
which can be written as $\mathbold x_k = \mathbold\phi_k^\top \mathbold d$ with $\mathbold{\phi}_k = [\mathbold x_{k-1}, \mathbold x_{k-2}, \cdots , \mathbold x_{k-m}]^\top$ and $\mathbold{d} = [d_0, \ldots, d_{m-1}]^\intercal \in \R^m$.
We can then apply recursive least squares to identify $\mathbold d$, along the same lines as~\eqref{eq:RLS_basic}.

We are now ready to show how to incorporate the projections into this identification scheme. We do so by changing the data we feed to the RLS. We continue to use the recurrence introduced in~\eqref{eq:p2_RLS_rec} to identify the internal model, however the input data is now the projected outputs of either online gradient descent or the control-based algorithm. In this case, assuming we remain inside the feasible region the assumptions used to identify the system hold. In the case that we remain outside of the constrained solution space for too long, the algorithm will simply fallback to POGD as the identified internal model will not be sufficiently accurate.

\section{Numerical Results}\label{sec:numerics}
In this section we validate the algorithm design proposed in section~\ref{sec:algorithm} for the microgrid management problem described in section~\ref{sec:problem}, comparing it to alternative algorithms.

\subsection{Simulation set-up}
\begin{figure}[!ht]
    \centering
    \includegraphics[width=\linewidth]{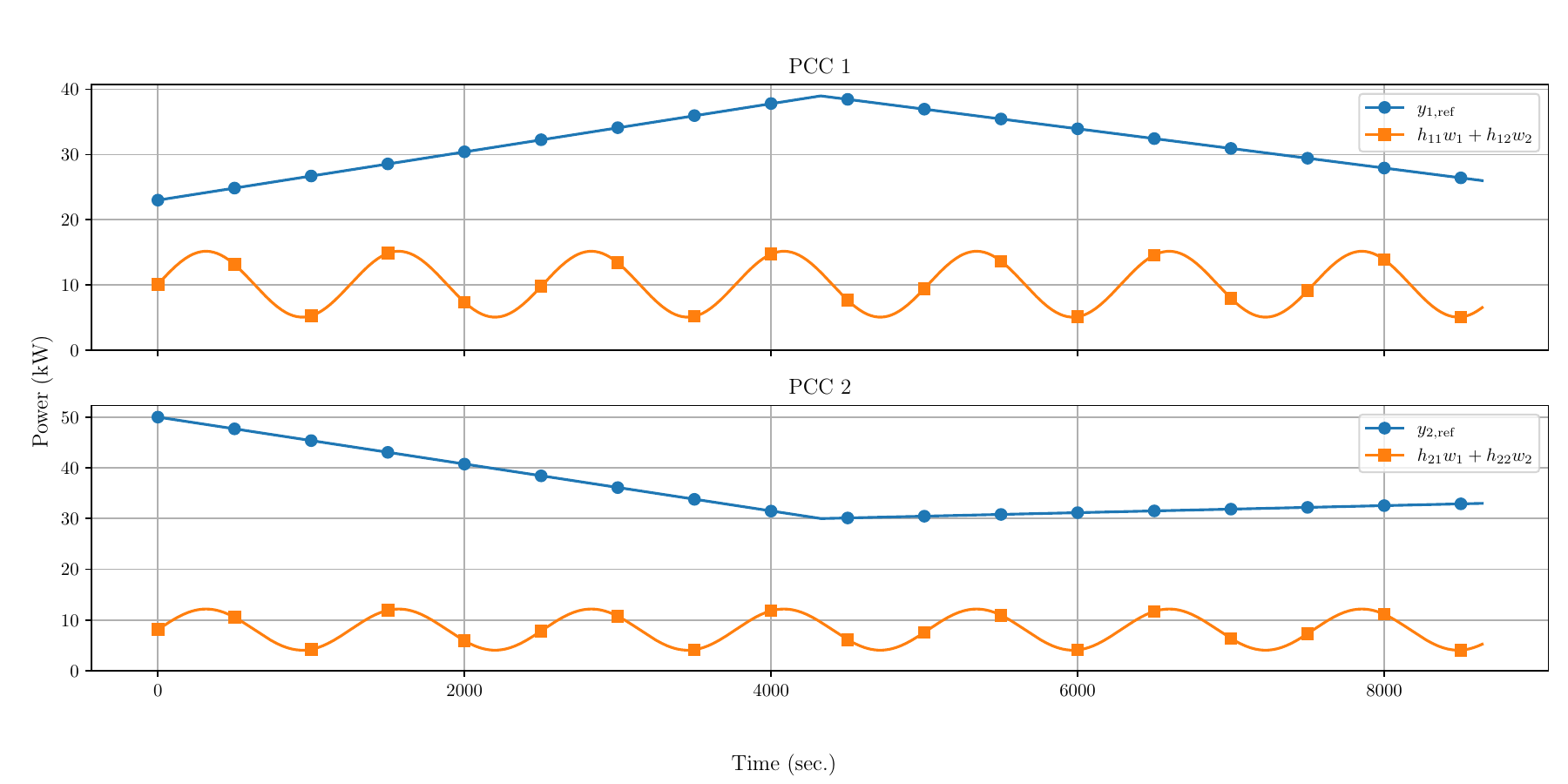}
    \caption{Reference signal $\boldsymbol y_{k,\text{ref}}$ and overall non-controllable load at PCCs}
    \label{fig:PCCs}
\end{figure}

As discussed in section~\ref{sec:problem}, the objective is to control some number of distributed energy resources within a microgrid, where there are a number of uncontrollable loads present, as well as some number of points of common coupling (PCCs).
More concretely, we consider $6$ DERs and $2$ PCCs, resulting in the following problem dimensions: $\boldsymbol x \in \mathbb{R}^6,\hspace{1ex} \boldsymbol w \in \mathbb R^2, \hspace{1ex} \boldsymbol J \in \mathbb{R}^{2\times 6},\hspace{1ex} \boldsymbol H \in \mathbb{R}^{2\times 2},\hspace{1ex} \boldsymbol y_{k,\text{ref}} \in \mathbb{R}^2 $.
Furthermore, as we consider the context of demand response in power systems, we consider the uncontrollable loads to be periodic (\textit{e.g.} each day has a similar user demand profile) \citep{Ospina}; in particular: $\boldsymbol w_k = \sin(0.005 k)$.
The reference signal, $\boldsymbol y_{k,\text{ref}}$, is generated as a portion of a triangular wave. Both the uncontrollable load and the reference signal at each of the two PCCs is shown in Figure \ref{fig:PCCs}.

The matrices $\boldsymbol H$ and $\boldsymbol J$ are randomly generated to ensure a realistic problem formulation, while guaranteeing that Assumption~\ref{as:cost-constraint} holds. In particular, $\boldsymbol H$ and $\boldsymbol J^\intercal \boldsymbol J$ are generated positive semi-definite with eigenvalues in $[1, 5]$.

These settings above then define the ``engineering" part of the cost in~\eqref{eq:true_objective}, and we need to define the user dissatisfaction part.
This function may represent, for example, the charging rate of an electric vehicle, or relative indoor temperature if the DER is an AC unit \citep{Ospina}.
Mathematically, we define the following model for $U_k(\boldsymbol x)$:
\begin{equation}
    U_k(\boldsymbol x) = \x^\intercal \boldsymbol{U}_{1} \x + \x^\intercal \boldsymbol{u}_{2,k} + u_{3,k}
\end{equation}
where the values of $\boldsymbol{U}_{1} \in \R^{6 \times 6}$, $\boldsymbol{u}_{2,k} \in \R^6$, $u_{3,k} \in \R$ are unknown and potentially time-varying. $\boldsymbol{U}_{1}$ is generated as a positive definite matrix with eigenvalues in $(0, 1]$. $\boldsymbol{u}_{2,k}$ and $u_{3,k}$ instead switch between two randomly generated values (taken from $\mathcal{N}(0,1)$) to model changes in the preferences of the DER owners. The first value is used in the first and third quarters of the simulation, the second value in the remaining times.
The resulting cost of~\eqref{eq:true_objective} is then a strongly convex quadratic cost, whose Hessian has eigenvalues in $[\lmin, \lmax] = [1,6]$.

To conclude, we need to define the time-varying constraint set $\mathcal{X}_k$.
In particular, we impose the following power limits at the DERs: $\mathcal{X}_k = [-10, 10] \times [-6, 6] \times [3, 13] \times [7, 17] \times [0, 28] \times [3, 32]$.

\subsection{Simulation results}
We are now ready to test the proposed algorithm in the set-up of the previous section.
In particular, we compare the proposed algorithm to the baseline (unstructured) projected online gradient descent (POGD) alternative, as well as to the projected control-based algorithm with anti-windup (PCBW) defined by~\eqref{eq:proj_CB}. The proposed algorithm, P-SIMBO, refers to the algorithm represented by Figure~\ref{fig:flowchart}, and introduced in section~\ref{sec:algorithm}.
\footnote{We remark that we tested projected SIMBO with the addition of anti-windup, but the results were not satisfying. This is an interesting direction of future research.}
We evaluate the algorithms' performance in terms of tracking error $\norm{\x_k - \x_{*,k}}$ where the optimal solutions $\x_{*,k}$ are computed via static optimization at each timestep using \texttt{cvxpy} \citep{diamond2016cvxpy}.
The simulations were implemented using \texttt{tvopt} \citep{tvopt}.

\begin{figure}[ht]
    \centering
    \includegraphics[width=\linewidth]{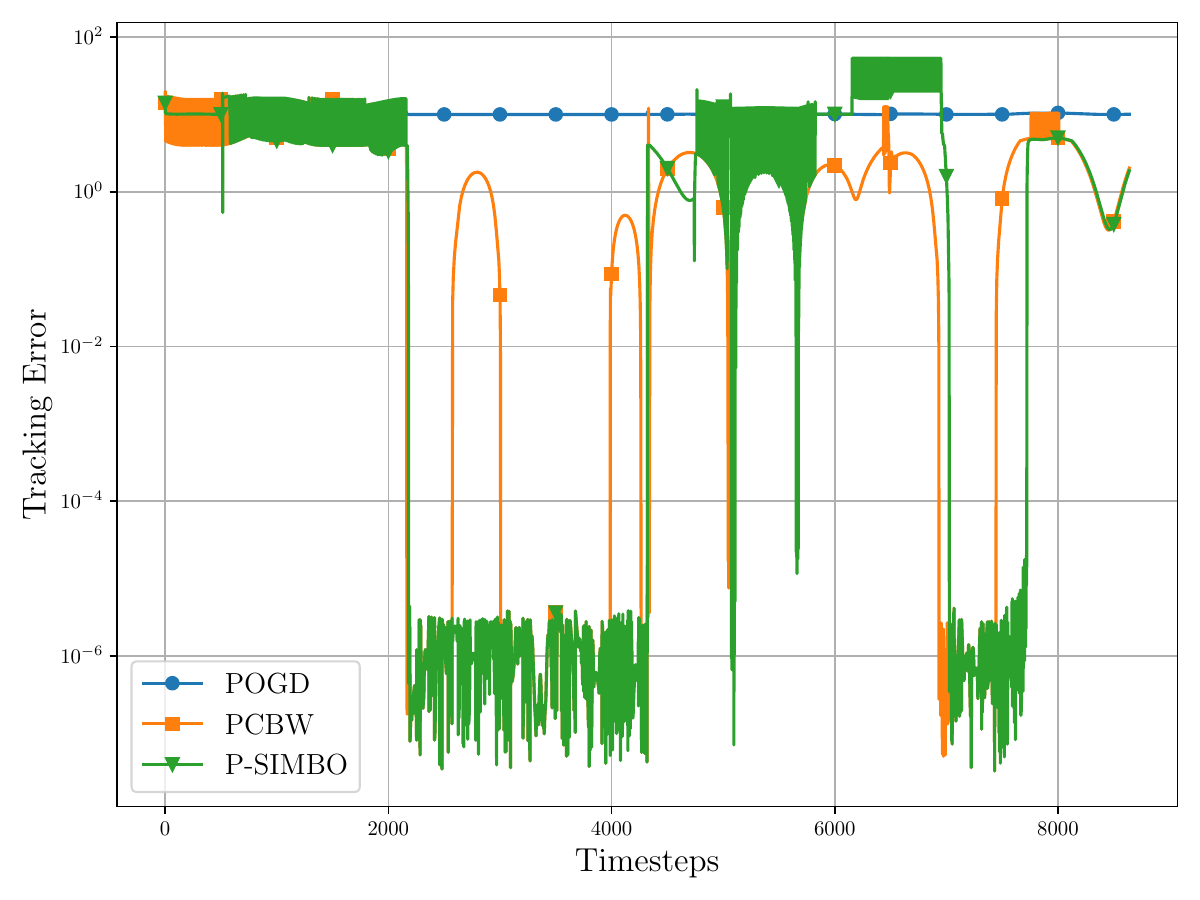}
    \caption{Tracking error comparison}
    \label{fig:tracking_error}
\end{figure}

We start by presenting in Figure~\ref{fig:tracking_error} the evolution of the tracking error across the simulation.
For the control-based algorithm, which does not employ identification, we use the internal model of a sine squared, given as:
\begin{align}
    \mathcal Z\{\sin^{2}(\omega_{0} n)\} = \frac{\sin^{2}(\omega_{0})\, z (z+1)}{(z-1)\bigl(z^{2}-2\cos(2\omega_{0})z+1\bigr)},
\end{align}
where $\omega_0=10.$ based on testing.
As we can see from Figure~\ref{fig:tracking_error}, initially both control-based algorithms perform similarly to the baseline of POGD. This is due to the fact that both control-based algorithms are undergoing a transient where performance is dominated by the sub-optimality of the initial conditions (see similar behaviors in \citep{Bastianello2024,casti2023}).
Notice also that P-SIMBO coincides with POGD for the first $\sim 500$ iterations (phase 1 in Figure~\ref{fig:flowchart}).
After $\sim 2000$ timesteps the user dissatisfaction changes, and both control-based algorithms significantly outperform POGD, both achieving tracking errors close to $10^{-6}$. Notice in particular that the tracking error of P-SIMBO consistently stays small.
Subsequent changes in the user dissatisfaction cost degrade the performance after $\sim 4000$ timesteps, but for the most part both control-based algorithms still outperform the POGD baseline. Transients where P-SIMBO yields worse performance are due to the internal model not being precisely identified yet after a change in the cost function.

\begin{figure}[ht]
    \centering
    \includegraphics[width=\linewidth]{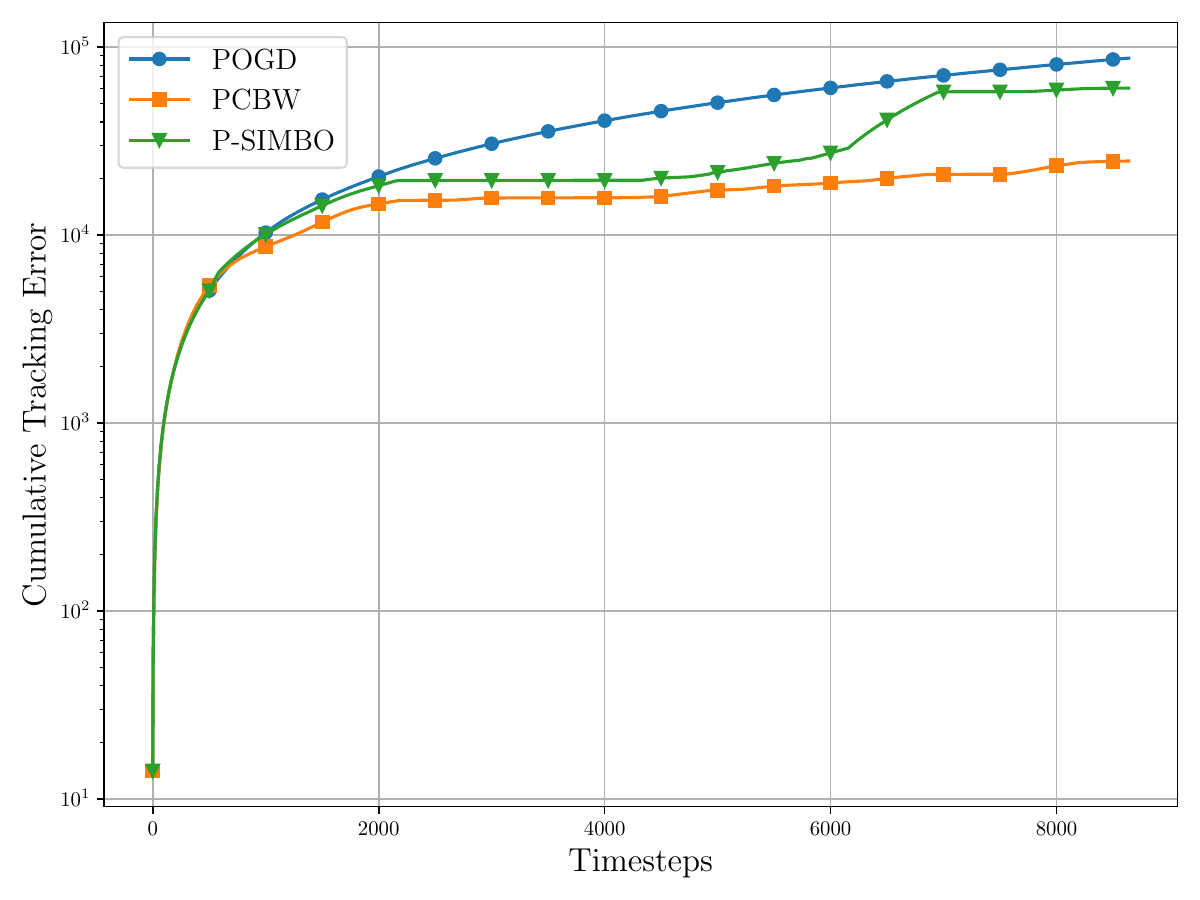}
    \caption{Cumulative tracking error comparison}
    \label{fig:cumulative_tracking_error}
\end{figure}

Clearly a realistic simulation such as this deviates from the theoretical assumptions made in section~\ref{subsec:algorithm}, yielding a more complex behavior than \textit{e.g.} simulations in \citep{Bastianello2024}.
Nonetheless, inspecting the long-term performance of the algorithms still allows us to highlight the improvement that control-based algorithms provide over the unstructured POGD.
In particular, Figure~\ref{fig:cumulative_tracking_error} depicts the cumulative tracking error across the simulation.
We can see that, as observed before, there is an initial transient with similar performance for all algorithms.
Then, as time progresses, the control-based algorithms outperform POGD and continue to do so in the long term. This demonstrates that a tailored online algorithm design has better performance than an unstructured algorithm, even in realistic scenarios which do not satisfy modeling assumptions precisely.
We remark that between timesteps $2000$ and $6000$, P-SIMBO and PCBW have largely comparable performance, while P-SIMBO's performance degrades afterwards, but still outperforming POGD. This later degradation of the performance, as discussed above, is due to changes in the user dissatisfaction cost which in turn require re-identification of the model. And while the model is not yet correctly identified, the performance is worse.

One important point to notice, however, is that P-SIMBO is entirely reliant on data about the online problem that is collected in an online fashion. On the other hand, PCBW owes its better performance on a more exact internal model that was selected offline from ``historical" data on the problem. But in many practical scenarios, historical data to identify a model might not be available, or the model might change over time. For this reason, identifying the internal model online as in P-SIMBO requires less information on the problem, and ensures adaptability to rapidly changing contexts, such as those which arise in microgrids with a high penetration of renewables.
In this sense, PCBW should be seen as an optimal baseline, which the adaptive algorithm would match if the microgrid conditions were to remain constant. The strength of P-SIMBO then lies in its adaptability to the problem.

\section{Conclusions and Future Work}\label{sec:conclusion}
In this paper we develop a novel online algorithm for problems with a (time-varying) convex constraint set, motivated by the management of microgrids with renewable sources.
The algorithm is characterized by a control-based design that applies the internal model principle, and a system identification routine tasked with identifying such internal model. In addition, in order to ensure the constraints are verified, we integrate a projection onto the constraint set.
We showcase promising numerical results for the microgrid use case, highlighting in particular the enhanced adaptability of the proposed algorithm to changes in the internal model.
Future work will address the use of more complex internal models (\textit{e.g.} non-linear or stochastic) to design the algorithm, with the goal of improving performance for realistic problems such as the management of a microgrid.

\bibliography{ifacconf}
\end{document}